\documentclass[reqno]{amsart}
\usepackage{amsmath, amsthm, amsfonts, amsopn, graphics}

    \theoremstyle{plain}
    \newtheorem{Thm}{Theorem}[section]
    \newtheorem{Prop}[Thm]{Proposition}
    
    \newtheorem*{Lemma*}{Lemma}
    \newtheorem{Cor}[Thm]{Corollary}

    \theoremstyle{definition}
    \newtheorem{Def}[Thm]{Definition}
    \newtheorem*{Def*}{Definition}
    \newtheorem{Example}[Thm]{Example}

    \theoremstyle{remark}
    
    \newtheorem*{Remark*}{Remark}

    \numberwithin{equation}{section}

    \newcommand{\field}[1]{\mathbb{#1}}
    \newcommand{\N}{\field{N}}
    \newcommand{\Z}{\field{Z}}

    \newcommand{\C}{\field{C}}
    \DeclareMathSymbol{\fieldk}{\mathalpha}{AMSb}{"7C} 
    \newcommand{\abs}[1]{\lvert#1\rvert}
    
    \newcommand{\norm}[1]{\lVert#1\rVert}
    \newcommand{\inner}[2]{\left\langle#1,#2\right\rangle}
    \DeclareMathOperator{\Tr}{Tr}
    \DeclareMathOperator{\Det}{Det}
    
    \DeclareMathSymbol{\normal}{\mathord}{AMSa}{"43}
    \DeclareMathOperator{\Log}{Log}
    \DeclareMathOperator{\Exp}{Exp}

    \newcommand{\bs}{\backslash}

    \newcommand{\ep}{\varepsilon}
    \newcommand{\la}{\lambda}
    \newcommand{\Laplace}{\Delta}
    
    \newcommand{\Del}{\Delta}
    \newcommand{\del}{\delta}

    \newcommand{\define}[1]{\emph{#1}}      

    \newcommand{\ltwo}{l^{2}}
    \newcommand{\Ltwo}{L^{2}}
    
    \newcommand{\VN}[1]{\mathcal{#1}}
    
    \newcommand{\Luck}{L\"uck}

    \newcommand{\len}{\ell}   
    \newcommand{\adj}{\delta}

    \newcommand{\grouptower}{{\pi = \pi_{1} \supset \pi_{2} \supset \dotsb}}
    \newcommand{\myint}{\int_{-k^-}^{k}}
\begin{document}

%
%
\title{Convergence of zeta functions of graphs}

\author{Bryan Clair \and Shahriar Mokhtari-Sharghi}
\address{CUNY Graduate Center,
        365 5th Avenue,
        New York, NY 10016}
\email{bclair@gc.cuny.edu}
\address{Department of Mathematics,
        Long Island University,
        Brooklyn Campus, 1 University Plaza,
        Brooklyn, NY 11201}
\email{mokhtari@liu.edu}
\date{\today}


\begin{abstract}
The $\Ltwo$-zeta function of an infinite graph $Y$ (defined previously
in a ball around  zero) has an analytic
extension.  For a tower of finite graphs
covered by $Y$, the normalized zeta functions of the finite graphs
converge to the $\Ltwo$-zeta function of $Y$.
\end{abstract}

\maketitle

%
%
\section*{Introduction}
Associated to any finite graph $X$ there is a zeta function $Z(X,u)$,
$u \in \C$.  It is defined as an
infinite product but shown (in various different cases) by Ihara,
Hashimoto, and Bass  \cite{Ih1, Ha1, Ba_zeta} to be a
polynomial. Indeed the rationality formula for a $q+1$ regular $X$ states that:
\begin{equation}\label{eq:rational}
Z(X,u)= (1-u^2)^{-\chi (X)} \Det(I - \del u + q u^2).
\end{equation}
Here $\del$ is the adjacency operator of $X$.

In \cite{cms:zeta},
an $\Ltwo$-zeta function is defined for  noncompact
graphs with symmetries, using the 
machinery of von Neumann algebras.  A rationality formula similar
to \eqref{eq:rational} expresses the relationship  between the zeta
function and 
the von Neumann determinant of a Laplace operator.
The results of this paper
focus on an especial case. 
Let $Y$ be an infinite graph which covers a finite graph $B=\pi\bs Y$.
The $\Ltwo$-zeta function  $Z_\pi(Y,u)$ is 
defined in \cite{cms:zeta} only in a small neighborhood of zero.  The first
result of this paper is to extend the $\Ltwo$-zeta function to the
interior of $C=\{u\in \C : \abs u=q^{-1/2} \} \cup [ -1 , - \frac 1
q]\cup [\frac 1 q , 1]$. 

In the second part of the paper, we consider a tower of finite
graphs $B_{i}$ covered by $Y$.  Put $N_i=\abs {B_i}/\abs{B}$. In
Theorem~\ref{thm:conv} we show 
that the zeta functions for the $B_{i}$, renormalized by taking
$N_i^{\text{th}}$ roots, converge to the $\Ltwo$-zeta function for $Y$.
The argument is inspired by, and uses, work of
\Luck~\cite{luck:resfin}.

In the first section we recall the 
definitions of the zeta functions of finite and infinite
graphs. One of the main results of this paper is
Theorem~\ref{extension} in this section.  
In the second section of this paper we prove the convergence theorem
and exhibit interesting examples. Theorem \ref{deit} 
generalizes work of Deitmar  \cite{deitmar}.

\section{Zeta functions}
In this section we recall the definition of the zeta function and
related material. We first recall the definition of the zeta function
for finite graphs.

\subsection{Finite Graphs}
For a graph $X$, let $VX$ and $EX$ denote the sets of vertices and edges,
respectively,  of $X$.
If each vertex has the same degree then $X$ is \define{regular}. 
\begin{Def}
Let $X$ be a finite graph. A closed path in $X$ is \define{primitive}
if it is  not a nontrivial  power of 
another path  inside the fundamental group of $X$. Let $P$ be
the set of free homoptopy 
classes of primitive
closed paths of $X$. Then the zeta function of $X$ is
\[
Z(X,u)=\prod_{\gamma\in P} (1 - u^{\len(\gamma)}),
\]
where $\len(\gamma)$ is the minimum length of paths in the class of $\gamma$. 
\end{Def}
Let $\del$ be the adjacency operator of the graph $X$ acting on
$\ltwo (VX)$. For $f\in \ltwo (VX)$ let $Qf(x)=q(x)f(x)$ where $q(x)+1$
is the degree of the vertex $x$.
Put $\Del(X,u)=I-\del u + Q u^2$.
The Ihara rationality formula says that 
$Z(X,u)$ is  a polynomial:
\begin{equation}
Z(X,u) = (1-u^2)^{-\chi(X)} \Det(\Del(X,u)).
\end{equation}
The zeta function satisfies the following  functional equation:
\cite[Corollary 3.10]{Ba_zeta}
\begin{equation}\label{eq:functional}
Z(X,(qu)^{-1}) =
  \left(\frac{1-u^{2}}{q^{2}u^{2}-1}\right)^{\chi(X)}
  q^{v-2e}u^{-2e}Z(X,u)
\end{equation}
where $e=\abs{EX}$ and $v=\abs{VX}$.

For more details and examples see \cite{starkter}.

The following proposition is well known. A proof can be found in
\cite[Page 59]{lubotzky}.
\begin{Prop}
The zeros of the zeta function for any finite $q+1$ regular graph
lie on the set $C$ where
\begin{equation*}
C=\{u\in \C : \abs u=q^{-1/2} \} \cup [ -1 , - \frac 1 q]\cup [\frac 1
q , 1].
\end{equation*}
\end{Prop}
\begin{figure}[h]   
   \begin{center}
    \input{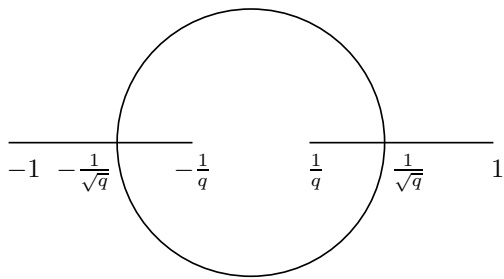}
    \caption{The set $C$}
    \label{pic:ci}
   \end{center}
\end{figure}
Let $\Omega$ be  the interior of $C$. Then:
\begin{Prop}\label{omega}
If $u \in \Omega, \lambda \in [-(q+1), q+1]$ then  $1-\lambda u +
qu^2\in \C- (-\infty,0]$.
\end{Prop}
\begin{proof}
Let $u=a +b i$. For $1-\lambda u + q u^2$ to be real we must have $1-\lambda u+
qu^2=1-\lambda \bar u +q {\bar u}^2.$  This implies
that either $b=0$ or $a=\frac{\la}{2q}$.
If $a=\frac{\la}{2q}$ then $1-\la u+q u^2= 1-q(a^2+b^2)$ which is
negative or zero if $u$ is on or outside of the circle $\abs
u=q^{-1/2}$.
In the case that $b=0$ put $f(a):=1-\la a+q a^2=1- \la u +q u^2$.
If $f$ has no real root then $f$ is always
positive. Otherwise $\abs \la \ge 2\sqrt q.$ Then $f(-\frac 1 q)$
and $f( \frac 1 q)$ are non-negative, and  $f'$ is not
 zero in
$(-\frac 1 q , \frac 1 q)$. Therefore $f$ will be positive on $(-\frac 1 q ,
\frac 1 q)$.
\end{proof}
\begin{Cor}\label{nthroot}
For a $q+1$ regular finite graph X, the polynomial $\Det \Del(X,u)$ has
an analytic
$N^{\text{th}}$  root on
$\Omega$ for all $N\in \N$.
\end{Cor}
\begin{proof}
We know $\Det \Del(X,u) =\Pi_{\la}(1-\la u+qu^2)$ where $\la$ varies over
eigenvalues of the adjacency operator of the graph. Then $
\Pi_{\la}(1-\la u+qu^2)^{\frac 1 N}$ is an analytic $N^{\text{th}}$
root for $\Det \Del(X,u)$.
\end{proof}

\subsection{Infinite Graphs}
For $\pi$ a countable discrete group, the \define{von Neumann algebra}
of $\pi$ is the algebra $\VN{N}(\pi)$ of bounded $\pi$-equivariant operators
from $\ltwo(\pi)$ to $\ltwo(\pi)$.  The \define{von Neumann trace} of
an element $f \in \VN{N}(\pi)$ is defined by
\[ \Tr_{\pi}f = \inner{f(e)}{e} \]
for $e \in \pi$ the unit element.  For $H =
\oplus_{i=1}^{n}\ltwo(\pi)$ and a bounded $\pi$-equivariant
operator $f:H \to H$,
define
\[ \Tr_{\pi}f = \sum_{i=1}^{n}\Tr_{\pi}f_{ii}. \]
The trace as defined is independent of the decomposition of $H$.
The von Neumann trace extends to bounded $\pi$-equivariant operators
on Hilbert $\VN{N}(\pi)$-modules, but they will not be needed.

Now let $Y = (VY,EY)$ be an infinite graph.
Suppose the group $\pi$ acts freely on $Y$ with
finite quotient $B$.

Let $P$ denote the set of
free homotopy classes of primitive closed paths in $Y$.
For $\gamma\in P$, $\len(\gamma)$ is the length of the shortest
representative of $\gamma$.
The group $\pi_{\gamma}$ is the stabilizer of $\gamma$ under the
action of $\pi$.
The \define{$\Ltwo$-zeta function} of $Y$ is defined in
\cite{cms:zeta} as the infinite product
\[
   Z_{\pi}(Y,u) = \prod_{\gamma \in \pi\bs P}
                   \left(1-u^{\len(\gamma)}\right)
                   ^{\frac{1}{\abs{\pi_{\gamma}}}}.
\]

The adjacency operator $\adj$ and  Laplace operator
$\Del(Y,u) = I - \adj u + Qu^{2}$ are $\pi$-equivariant operators
on $\ltwo(VY)$.
Choosing lifts of vertices of $B$ yields a decomposition
$\ltwo(VY) = \oplus\ltwo(\pi)$.
Then from \cite[Theorem 0.3]{cms:zeta},
\begin{equation}\label{eq:zetadef}
   Z_{\pi}(Y,u) = (1-u^{2})^{-\chi^{(2)}(Y)}\Det_{\pi}\Del(Y,u).
\end{equation}
In this formula, $\chi^{(2)}(Y)$ is the $\Ltwo$-Euler characteristic
of $Y$, which in our setting is simply equal to $\chi(B)$.
The determinant $\Det_{\pi}\Laplace(Y,u)$ is defined via formal power
series as
$(\Exp \circ Tr_{\pi} \circ \Log) \Laplace(Y,u)$ and converges for
small $u \in \C$.
More precisely, if $Y$ is $q+1$ regular then the radius of convergence
of $Z_{\pi}(Y,u)$ is greater than or equal to $\frac 1 q.$

\begin{Thm}\label{extension}
Let $Y$ be a $q+1$ regular graph. Then $Z_{\pi}(Y,u)$ has a
holomorphic extension to $\Omega.$
\end{Thm}
\begin{proof}
By \eqref{eq:zetadef} it is enough to show that $\Det_\pi\Laplace(Y,u)$
has a holomorphic extension on $\Omega$.
Let $g_u(\lambda)=\log(1 - \lambda u + q u^2)$.
Here and in the rest of the paper $\log$ is the principal branch of
the logarithm, defined and analytic on $\C-(-\infty,0]$.
Fix $u\in\Omega$.  Then using Proposition~\ref{omega},
there exists an open set $V_u\supset [-(q+1),q+1]$ on which
$g_u$ is analytic.
Since $\adj$ is self-adjoint
and $\norm \adj \le q+1$, the spectrum $\sigma(\adj)\subset [-(q+1),q+1]$.
By the spectral theorem for self-adjoint operators we can write:
\begin{equation}
   \adj = \myint \lambda\ dE(\lambda).
\end{equation}
Now
\begin{equation}
   g_u(\adj) = \myint \log (1- \lambda  u + q u^2)\ dE(\lambda)
\end{equation}
is well defined, and
\begin{equation}\label{eq:tr_gu}
   \Tr_\pi g_u(\adj) = \myint
                       \log (1- \lambda  u + q u^2)\ d(\Tr_\pi E(\lambda))
\end{equation}
is a holomorphic function of $u$ on $\Omega$.
Now for small $u$,
\begin{equation}
   \begin{align}
      \Det_\pi(\Laplace(Y,u)) &= \Exp\Tr_\pi\Log(\Laplace(Y,u))\\
                            &= \Exp\Tr_\pi\Log(I - \adj u + q u^2)\\
                            &= \Exp\Tr_\pi g_u(\adj).
   \end{align}
\end{equation}
\end{proof}

\begin{Remark*}
Using the functional equation \eqref{eq:functional} it is possible to
extend the $L^2$-zeta 
function to the exterior of $C$.  
Defining the zeta function on $C$ itself presents some problems. In
examples, the function on the interior of $C$ and the exterior of $C$
do not match continuously on $C$. 
However the absolute value of the resulting function is the
Fuglede-Kadison determinant, and may extend continuously to $C$.
\end{Remark*}

\section{Convergence of Zeta Functions for Towers of Graphs}
In this section we prove that the zeta functions for a tower of finite
graphs, suitably renormalized, converge to the $\Ltwo$-zeta function
for an infinite covering graph.  The argument uses an idea from
\cite{luck:resfin}.
\subsection{The Convergence Theorem}
\begin{Thm}\label{thm:conv}
Let $Y$ be a $q+1$ regular graph. Suppose
the group $\pi$ acts freely on $Y$, and that
$B=\pi\bs Y$ is a finite graph.
Suppose $\grouptower$ is a tower of finite index
normal subgroups and $\bigcap\pi_i =\{e\}$.
Let $[\pi:\pi_i]=N_i$ and $B_i=\pi_i\bs Y$. Then for $u \in \Omega$,
we have
\begin{equation}
  \lim_{i\to \infty} Z(B_i,u)^{\frac 1 {N_i}} = Z_\pi(Y,u).
\end{equation}
The convergence is uniform on compact subsets of $\Omega$.
\end{Thm}
\begin{Remark*}
Notice 
\begin{equation}
\begin{align}
Z(B_i,u)&= (1-u^2)^{-\chi(B_i)}\Det\Del(B_i,u)\\
         &= (1-u^2)^{-N_i \chi(B)}\Det\Del(B_i,u).
\end{align}
\end{equation}
So we let 
\begin{equation}\label{limeq}
Z(B_i,u)^{1/N_i}= (1-u^2)^{-\chi(B)} (\Det\Del(B_i,u)) ^{1/N_i}
\end{equation}
The $N^{\text{th}}$ root in \eqref{limeq} is taken in the
sense of Corollary~\ref{nthroot}.
\end{Remark*}
\begin{proof}[Proof of Theorem \ref{thm:conv}]
From the remark we need to show that 
\begin{equation}
\lim_{i\to\infty}(\Det\Del(B_i,u)) ^{1/N_i}=
\Det\Del(B,u).
\end{equation}
Let $F_i(\la)=\frac{1}{N_i}{\abs {\{\text{$\mu$ eigenvalue of $\del_i$,
$\mu\le\la$} \} }}$.
Let $F(\lambda)=\Tr_\pi E(\la)$, where $\{E(\la)\}_\la$ is the
spectral decomposition of $\del$ acting on $\ltwo (Y).$
We now set
\[ \begin{array}{rclrcl}
      \overline{F}(\lambda) & = &
      		\limsup_{i \rightarrow \infty} F_{i}(\lambda); &
      \underline{F}(\lambda) & = &
      		\liminf_{i \rightarrow \infty} F_{i}(\lambda); \\
      \overline{F}^{+}(\lambda) & = &
      		\lim_{\ep \rightarrow 0^{+}}\overline{F}(\lambda+\ep);&
      \underline{F}^{+}(\lambda) & = &
      		\lim_{\ep \rightarrow 0^{+}}\underline{F}(\lambda+\ep).\\
      
\end{array} \]
Then from \cite[Theorem 2.3.1]{luck:resfin},
for all $\lambda \in [-k,k]$,
	\[ F(\lambda) = \overline{F}^{+}(\lambda)
						=
\underline{F}^{+}(\lambda)  \] 

We know
\begin{equation}\log\Det \Del(Y_i, u)=\myint\log
(1-u\la+qu^2)dF_i(\la).
\end{equation} By \eqref{eq:tr_gu}
\begin{equation}
  \log\Det_\pi\Del(Y,u) = \myint\log(1-u\la+qu^2)dF(\la).
\end{equation}
If $K \subset \Omega$ is compact then $\log(1-u\la+qu^2)$ is
bounded uniformly for $u \in K$ and $\la$ in
an open interval containing $[-k, k].$
Now integration by parts shows that indeed 
\begin{equation}\myint\log
(1-u\la+qu^2)dF_i(\la)\to \myint\log(1-u\la+qu^2)dF(\la).
\end{equation}
as  $i\to\infty$. 
\end{proof}

\subsection{Examples}
\begin{Example}
Let $Y$ be the line, with vertices $\Z$ and edges connecting $n$ to
$n+1$.  The group $\pi = \Z$ acts on $Y$ with quotient having one
vertex and one edge.  Let $\pi_{n} = n\Z \subset \Z$, so
$B_{n} = \pi_{n}\bs Y$ is an $n$-cycle.  As $B_{n}$ has only two
primitive loops and $Y$ has none, we have
$Z(B_{n},u) = (1-u^{n})^{2}$ and $Z_{\pi}(Y,u) = 1$.
The graphs are 2-regular, so $\Omega$ is the (open) unit disk. 
For $u \in \Omega$,
\[
    \lim_{n\to\infty}(1-u^{n})^{2/n} = 1.
\]
Notice that the functional equation \eqref{eq:functional}
gives $\Z_\pi(Y,u)=u^2$ outside of the unit disk.
\end{Example}

\begin{Example}
Let $Y$ be the $q+1$ regular tree and $B = \pi\bs Y$ finite.
Since $Y$ has no closed loops, the $\Ltwo$-zeta function of $Y$ is the
constant function 1.
So if $B_{n} = \pi_{n}\bs Y$ is a tower
of finite graphs covering $B$,
\[
   \lim_{n\to\infty}Z(B_{n},u)^{1/[\pi:\pi_{n}]} = 1
\]
for $u \in \Omega$.
This result is contained in \cite{deitmar} for $\abs{u}$ small.
\end{Example}

\begin{Thm}\label{deit}
Let $Y$ and $B$ be as in the previous example.  We have:
\[
Z(B,u)=\frac{\Det\Del(B,u)}{\Det_\pi\Del(Y,u)}
\]
for $u\in\Omega$.
\end{Thm}
For small $u$, the above theorem is the main result of
\cite{deitmar}.
\begin{proof}
From \cite[Theorem 0.3]{cms:zeta} we know that
\[
Z_\pi(Y,u)=(1-u^2)^{-\chi^{(2)}(Y)}\Det_\pi\Del(Y,u).
\]
As in the previous example, $Z_\pi(Y,u)=1$.  Now we have
\[
  Z(B,u) = \frac{Z(B,u)}{Z_\pi(Y,u)}
         = \frac{(1-u^2)^{-\chi(B)}\Det\Del(B,u)}
                {(1-u^2)^{-\chi^{(2)}(Y)}\Det_\pi\Del(Y,u)}.
\]
As $\chi(B)=\chi^{(2)}(Y)$, the theorem follows.
\end{proof}

%
%
\bibliographystyle{plain}

\end{document}